\documentclass[12pt,a4paper]{article}

\usepackage{graphicx}

\usepackage{amsmath}

\usepackage{amssymb}
\usepackage{hyperref}
\usepackage{authblk}

\usepackage{amsthm}
\usepackage{geometry}
\usepackage{float}

\usepackage{setspace}
\spacing{1.5}
\usepackage{fancyhdr}

\usepackage{color} 

\graphicspath{{pics/}}

\newcommand{\al}{\alpha}

\newcommand{\vphi}{\varphi}

\newcommand{\be}{\beta}

\newcommand{\ga}{\gamma}

\newcommand{\de}{\delta}

\newcommand{\om}{\omega}

\newcommand{\na}{\nabla}

\newcommand{\NA}{\nabla}

\newcommand{\bs}{\boldsymbol}

\newcommand{\ra}{\rightarrow}

\newcommand{\lra}{\longrightarrow}

\newcommand{\Ra}{\Rightarrow}

\newcommand{\xra}{\xrightarrow}

\newcommand{\xlra}{\xlongrightarrow}

\newcommand{\rgl}{\rangle}

\newcommand{\lgl}{\langle}

\newcommand{\dash}{\textrm{-}}

\newcommand{\ot}{\otimes}

\newcommand{\bpf}{\begin{proof}}

\newcommand{\epf}{\end{proof}}

\newcommand{\bthm}{\begin{thm}}

\newcommand{\ethm}{\end{thm}}

\newcommand{\bprop}{\begin{prop}}

\newcommand{\eprop}{\end{prop}}

\newcommand{\bcor}{\begin{cor}}

\newcommand{\ecor}{\end{cor}}

\newcommand{\blem}{\begin{lem}}

\newcommand{\elem}{\end{lem}}

\newcommand{\bdefn}{\begin{defn}}

\newcommand{\edefn}{\end{defn}}

\newcommand{\bexmp}{\begin{exmp}}

\newcommand{\eexmp}{\end{exmp}}

\newcommand{\brem}{\begin{rem}}

\newcommand{\erem}{\end{rem}}

\newcommand{\bdia}{\begin{displaymath}\xymatrix}

\newcommand{\edia}{\end{displaymath}}

\newcommand{\beq}{\begin{equation*}\begin{aligned}}

\newcommand{\eeq}{\end{aligned}\end{equation*}}

\newcommand{\bref}{\textbf{Ref}}

\newcommand{\intg}{\mathbb{Z}}

\newcommand{\real}{\mathbb{R}}

\newcommand{\comp}{\mathbb{C}}

\newcommand{\quot}{\mathbb{H}}

\newcommand{\afv}{\mathbb{A}}

\newcommand{\prv}{\mathbb{P}}

\newcommand{\mco}{\mathcal{O}}

\newcommand{\mcc}{\mathcal{C}}

\newcommand{\mcf}{\mathcal{F}}

\newcommand{\mcg}{\mathcal{G}}

\newcommand{\mcs}{\mathcal{S}}

\newcommand{\cp}{\mathbb{CP}}

\newcommand{\mfo}{\mathfrak{O}}

\newcommand{\mfg}{\mathfrak{g}}

\newcommand{\msa}{\mathscr{A}}

\newcommand{\msr}{\mathscr{R}}

\newcommand{\msg}{\mathscr{G}}

\newcommand{\msd}{\mathscr{D}}

\newcommand{\itbf}{\item\textbf}

\newcommand{\seqa}{a_1,...,a_}

\newcommand{\seqx}{x_1,...,x_}

\newcommand{\seqy}{y_1,...,y_}

\newcommand{\seqf}{f_1,...,f_}

\newcommand{\cred}{\textcolor{red}}

\newcommand{\cblue}{\textcolor{blue}}

\newcommand{\mfa}{\mathfrak{a}}

\newcommand{\mfb}{\mathfrak{b}}

\newcommand{\mfm}{\mathfrak{m}}

\newcommand{\mfn}{\mathfrak{n}}

\newcommand{\mfp}{\mathfrak{p}}

\newcommand{\Af}{A_{(f)}}
\DeclareMathOperator{\si}{SI}
\DeclareMathOperator{\I}{i}
\DeclareMathOperator{\Int}{int}

\newtheorem{thm}{\textbf {Theorem}}[section]

\newtheorem{cor}[thm]{\textbf{Corollary}}

\newtheorem{prop}[thm]{\textbf{Proposition}}

\newtheorem{lem}[thm]{\textbf{Lemma}}

\newtheorem{conj}[thm]{Conjecture}

\theoremstyle{definition}

\newtheorem{defn}[thm]{\textbf{Definition}}

\newtheorem{exmp}[thm]{Example}

\newtheorem{notn}[thm]{Notation}

\theoremstyle{remark}

\newtheorem{rem}[thm]{Remark}

\def\cok{\operatorname{Coker}}

\newcommand{\txi}{\tilde{\xi}}

\newcommand{\bxi}{\bar{\xi}}

\newcommand{\bz}{\bar{z}}

\DeclareMathOperator{\tr}{trunk}

\title{Equivalence Relations Between Closed Curves On the Pair of Pants}

\author{Nithin Kavi}

\date{}

\def\allfiles{}
\begin{document}

\bibliographystyle{plain}
\maketitle

\begin{abstract}

We examine an equivalence relation between free homotopy classes of closed curves on the pair of pants known as k-equivalence, a generalization of a concept previously defined by Leininger. We prove that two classes of closed curves on the pair of pants that are k-equivalent must also be 1-equivalent and 2-equivalent. We also examine properties of 1-equivalence on the pair of pants in greater depth.

\end{abstract}

\tableofcontents


\ifx\allfiles\undefined

\documentclass[12pt,a4paper]{article}


\usepackage{graphicx}

\usepackage{amsmath}

\usepackage{amssymb}

\usepackage{amsthm}

\usepackage{geometry}

\DeclareMathOperator{\tr}{trunk}

\usepackage{fancyhdr}

\usepackage{color} 









\newcommand{\al}{\alpha}

\newcommand{\vphi}{\varphi}

\newcommand{\be}{\beta}

\newcommand{\ga}{\gamma}

\newcommand{\de}{\delta}

\newcommand{\om}{\omega}

\newcommand{\na}{\nabla}

\newcommand{\NA}{\nabla}

\newcommand{\bs}{\boldsymbol}

\newcommand{\ra}{\rightarrow}

\newcommand{\lra}{\longrightarrow}

\newcommand{\Ra}{\Rightarrow}

\newcommand{\xra}{\xrightarrow}

\newcommand{\xlra}{\xlongrightarrow}

\newcommand{\rgl}{\rangle}

\newcommand{\lgl}{\langle}

\newcommand{\dash}{\textrm{-}}

\newcommand{\ot}{\otimes}

\newcommand{\bpf}{\begin{proof}}

\newcommand{\epf}{\end{proof}}

\newcommand{\bthm}{\begin{thm}}

\newcommand{\ethm}{\end{thm}}

\newcommand{\bprop}{\begin{prop}}

\newcommand{\eprop}{\end{prop}}

\newcommand{\bcor}{\begin{cor}}

\newcommand{\ecor}{\end{cor}}

\newcommand{\blem}{\begin{lem}}

\newcommand{\elem}{\end{lem}}

\newcommand{\bdefn}{\begin{defn}}

\newcommand{\edefn}{\end{defn}}

\newcommand{\bexmp}{\begin{exmp}}

\newcommand{\eexmp}{\end{exmp}}

\newcommand{\brem}{\begin{rem}}

\newcommand{\erem}{\end{rem}}

\newcommand{\bdia}{\begin{displaymath}\xymatrix}

\newcommand{\edia}{\end{displaymath}}

\newcommand{\beq}{\begin{equation*}\begin{aligned}}

\newcommand{\eeq}{\end{aligned}\end{equation*}}

\newcommand{\bref}{\textbf{Ref}}

\newcommand{\intg}{\mathbb{Z}}

\newcommand{\real}{\mathbb{R}}

\newcommand{\comp}{\mathbb{C}}

\newcommand{\quot}{\mathbb{H}}

\newcommand{\afv}{\mathbb{A}}

\newcommand{\prv}{\mathbb{P}}

\newcommand{\mco}{\mathcal{O}}

\newcommand{\mcc}{\mathcal{C}}

\newcommand{\mcf}{\mathcal{F}}

\newcommand{\mcg}{\mathcal{G}}

\newcommand{\mcs}{\mathcal{S}}

\newcommand{\cp}{\mathbb{CP}}

\newcommand{\mfo}{\mathfrak{O}}

\newcommand{\mfg}{\mathfrak{g}}

\newcommand{\msa}{\mathscr{A}}

\newcommand{\msr}{\mathscr{R}}

\newcommand{\msg}{\mathscr{G}}

\newcommand{\msd}{\mathscr{D}}

\newcommand{\itbf}{\item\textbf}

\newcommand{\seqa}{a_1,...,a_}

\newcommand{\seqx}{x_1,...,x_}

\newcommand{\seqy}{y_1,...,y_}

\newcommand{\seqf}{f_1,...,f_}

\newcommand{\cred}{\textcolor{red}}

\newcommand{\cblue}{\textcolor{blue}}

\newcommand{\mfa}{\mathfrak{a}}

\newcommand{\mfb}{\mathfrak{b}}

\newcommand{\mfm}{\mathfrak{m}}

\newcommand{\mfn}{\mathfrak{n}}

\newcommand{\mfp}{\mathfrak{p}}

\newcommand{\Af}{A_{(f)}}


\newtheorem{thm}{\textbf {Theorem}}[section]

\newtheorem{cor}[thm]{\textbf{Corollary}}

\newtheorem{prop}[thm]{\textbf{Proposition}}

\newtheorem{lem}[thm]{\textbf{Lemma}}

\newtheorem{conj}[thm]{Conjecture}

\newtheorem{prob}[thm]{Problem}

\newtheorem{exer}[thm]{Exercise}

\newtheorem{quest}[thm]{Question}

\theoremstyle{definition}

\newtheorem{defn}[thm]{\textbf{Definition}}

\newtheorem{defns}[thm]{Definitions}

\newtheorem{exmp}[thm]{Example}

\newtheorem{exmps}[thm]{Examples}

\newtheorem{var}[thm]{Variant}

\newtheorem{vars}[thm]{Variants}

\newtheorem{con}[thm]{Construction}

\newtheorem{notn}[thm]{Notation}

\newtheorem{notns}[thm]{Notations}

\newtheorem{conv}[thm]{Convention}

\theoremstyle{remark}

\newtheorem{rem}[thm]{Remark}

\newtheorem{rems}[thm]{Remarks}

\newtheorem{warn}[thm]{Warning}

\newtheorem{sch}[thm]{Scholium}

\newtheorem{expl}[thm]{Explanations}

\newtheorem*{theorem}{\textbf{Theorem}}

\newtheorem*{corollary}{\textbf{Corollary}}

\newtheorem*{proposition}{\textbf{Proposition}}

\newtheorem*{lemma}{\textbf{Lemma}}

\newtheorem*{example}{\textbf{Example}}

\def\cok{\operatorname{Coker}}

\newcommand{\txi}{\tilde{\xi}}

\newcommand{\bxi}{\bar{\xi}}

\newcommand{\bz}{\bar{z}}



\begin{document}

\bibliographystyle{plain}

\maketitle

\else

\fi



\section{Introduction}

The pair of pants surface, or the triply punctured sphere, is an important topological surface. As the only surface with a finite number of deformation classes of closed curves of any fixed self intersection number, the pair of pants is the easiest surface on which to study properties of free homotopy classes of closed curves.

The concept of {\it simple-intersection equivalence} was first defined by Leininger in \cite{leininger}, and it groups deformation classes of curves into equivalence classes based on their intersections with classes of simple closed curves. On the pair of pants surface, however, all deformation classes of curves are simple-intersection equivalent. For this reason, we extend Leininger's definition to the concept of {\it k-equivalence}, which replaces classes of simple closed curves with deformation classes of closed curves with self-intersection number $k$ based on the following definition:

\textbf{Definition:} On a given surface $M$, two classes of curves $\al_1$ and $\al_2$ are defined to be {\it k-equivalent} if every class of curves $c$ on $M$ with self intersection number $k$ satisfies $i(\al_1, c) = i(\al_2, c).$

In this paper we study k-equivalence on the pair of pants, and prove the following result:

\textbf{Theorem:} On the pair of pants, if two classes of curves $\al_1$ and $\al_2$ are k-equivalent for some positive integer $k \geq 2,$ it follows that $\al_1$ and $\al_2$ are 2-equivalent and 1-equivalent.

This theorem also offers evidence for the following stronger conjecture:

\textbf{Conjecture 1:} On the pair of pants, if two classes of curves $\al_1$ and $\al_2$ are k-equivalent for some positive integer $k \geq 2,$ it follows that $\al_1$ and $\al_2$ are also $(k - 1)$ equivalent.

While analyzing 1-equivalence, we also found evidence for the following conjecture:

\textbf{Conjecture 2:} For any free homotopy class of curves $\al$ on the pair of pants, the ratio between the intersection numbers of $\al$ with one free homotopy class with self intersection number $1$ and another free homotopy class with self intersection number $1$ is always between $\frac{1}{2}$ and $2,$ inclusive.

The paper is organized as follows: In Section $2$, we introduce basic definitions and ideas related to closed curves on surfaces up to defining k-equivalence. Then, in Section $3,$ we see why the theorem above is true, and why conjectures $1$ and $2$ are likely true as well. We conclude with possible future directions of study and ideas on how to generalize the results presented here.

\section*{Acknowledgments}

I would like to thank my mentor Professor Moira Chas for offering guidance throughout this project and Professor Tao Li for providing feedback on my work. Additionally, I thank the Simons Summer Research Program, funded by the Simons Foundation, and Stony Brook University for giving me the opportunity to conduct this research.



\ifx\allfiles\undefined

\bibliography{Index}

\end{document}

\fi


\ifx\allfiles\undefined

\documentclass[12pt,a4paper]{article}


\usepackage{graphicx}

\usepackage{amsmath}

\usepackage{amssymb}

\usepackage{amsthm}

\usepackage{geometry}

\usepackage{fancyhdr}

\usepackage{color} 









\newcommand{\al}{\alpha}

\newcommand{\vphi}{\varphi}

\newcommand{\be}{\beta}

\newcommand{\ga}{\gamma}

\newcommand{\de}{\delta}

\newcommand{\om}{\omega}

\newcommand{\na}{\nabla}

\newcommand{\NA}{\nabla}

\newcommand{\bs}{\boldsymbol}

\newcommand{\ra}{\rightarrow}

\newcommand{\lra}{\longrightarrow}

\newcommand{\Ra}{\Rightarrow}

\newcommand{\xra}{\xrightarrow}

\newcommand{\xlra}{\xlongrightarrow}

\newcommand{\rgl}{\rangle}

\newcommand{\lgl}{\langle}

\newcommand{\dash}{\textrm{-}}

\newcommand{\ot}{\otimes}

\newcommand{\bpf}{\begin{proof}}

\newcommand{\epf}{\end{proof}}

\newcommand{\bthm}{\begin{thm}}

\newcommand{\ethm}{\end{thm}}

\newcommand{\bprop}{\begin{prop}}

\newcommand{\eprop}{\end{prop}}

\newcommand{\bcor}{\begin{cor}}

\newcommand{\ecor}{\end{cor}}

\newcommand{\blem}{\begin{lem}}

\newcommand{\elem}{\end{lem}}

\newcommand{\bdefn}{\begin{defn}}

\newcommand{\edefn}{\end{defn}}

\newcommand{\bexmp}{\begin{exmp}}

\newcommand{\eexmp}{\end{exmp}}

\newcommand{\brem}{\begin{rem}}

\newcommand{\erem}{\end{rem}}

\newcommand{\bdia}{\begin{displaymath}\xymatrix}

\newcommand{\edia}{\end{displaymath}}

\newcommand{\beq}{\begin{equation*}\begin{aligned}}

\newcommand{\eeq}{\end{aligned}\end{equation*}}

\newcommand{\bref}{\textbf{Ref}}

\newcommand{\intg}{\mathbb{Z}}

\newcommand{\real}{\mathbb{R}}

\newcommand{\comp}{\mathbb{C}}

\newcommand{\quot}{\mathbb{H}}

\DeclareMathOperator{\tr}{trunk}
\newcommand{\afv}{\mathbb{A}}

\newcommand{\prv}{\mathbb{P}}

\newcommand{\mco}{\mathcal{O}}

\newcommand{\mcc}{\mathcal{C}}

\newcommand{\mcf}{\mathcal{F}}

\newcommand{\mcg}{\mathcal{G}}

\newcommand{\mcs}{\mathcal{S}}

\newcommand{\cp}{\mathbb{CP}}

\newcommand{\mfo}{\mathfrak{O}}

\newcommand{\mfg}{\mathfrak{g}}

\newcommand{\msa}{\mathscr{A}}

\newcommand{\msr}{\mathscr{R}}

\newcommand{\msg}{\mathscr{G}}

\newcommand{\msd}{\mathscr{D}}

\newcommand{\itbf}{\item\textbf}

\newcommand{\seqa}{a_1,...,a_}

\newcommand{\seqx}{x_1,...,x_}

\newcommand{\seqy}{y_1,...,y_}

\newcommand{\seqf}{f_1,...,f_}

\newcommand{\cred}{\textcolor{red}}

\newcommand{\cblue}{\textcolor{blue}}

\newcommand{\mfa}{\mathfrak{a}}

\newcommand{\mfb}{\mathfrak{b}}

\newcommand{\mfm}{\mathfrak{m}}

\newcommand{\mfn}{\mathfrak{n}}

\newcommand{\mfp}{\mathfrak{p}}

\newcommand{\Af}{A_{(f)}}


\newtheorem{thm}{\textbf {Theorem}}[section]

\newtheorem{cor}[thm]{\textbf{Corollary}}

\newtheorem{prop}[thm]{\textbf{Proposition}}

\newtheorem{lem}[thm]{\textbf{Lemma}}

\newtheorem{conj}[thm]{Conjecture}

\newtheorem{prob}[thm]{Problem}

\newtheorem{exer}[thm]{Exercise}

\newtheorem{quest}[thm]{Question}

\theoremstyle{definition}

\newtheorem{defn}[thm]{\textbf{Definition}}

\newtheorem{defns}[thm]{Definitions}

\newtheorem{exmp}[thm]{Example}

\newtheorem{exmps}[thm]{Examples}

\newtheorem{var}[thm]{Variant}

\newtheorem{vars}[thm]{Variants}

\newtheorem{con}[thm]{Construction}

\newtheorem{notn}[thm]{Notation}

\newtheorem{notns}[thm]{Notations}

\newtheorem{conv}[thm]{Convention}

\theoremstyle{remark}

\newtheorem{rem}[thm]{Remark}

\newtheorem{rems}[thm]{Remarks}

\newtheorem{warn}[thm]{Warning}

\newtheorem{sch}[thm]{Scholium}

\newtheorem{expl}[thm]{Explanations}

\newtheorem*{theorem}{\textbf{Theorem}}

\newtheorem*{corollary}{\textbf{Corollary}}

\newtheorem*{proposition}{\textbf{Proposition}}

\newtheorem*{lemma}{\textbf{Lemma}}

\newtheorem*{example}{\textbf{Example}}

\def\cok{\operatorname{Coker}}

\newcommand{\txi}{\tilde{\xi}}

\newcommand{\bxi}{\bar{\xi}}

\newcommand{\bz}{\bar{z}}

\DeclareMathOperator{\tr}{tr}



\begin{document}

\bibliographystyle{plain}

\else

\fi

\section{Preliminaries}

We begin by giving basic definitions, much of which can be found in greater detail in \cite{zhang}.



\begin{figure}[H]
    \centering
    \includegraphics[scale = 0.2]{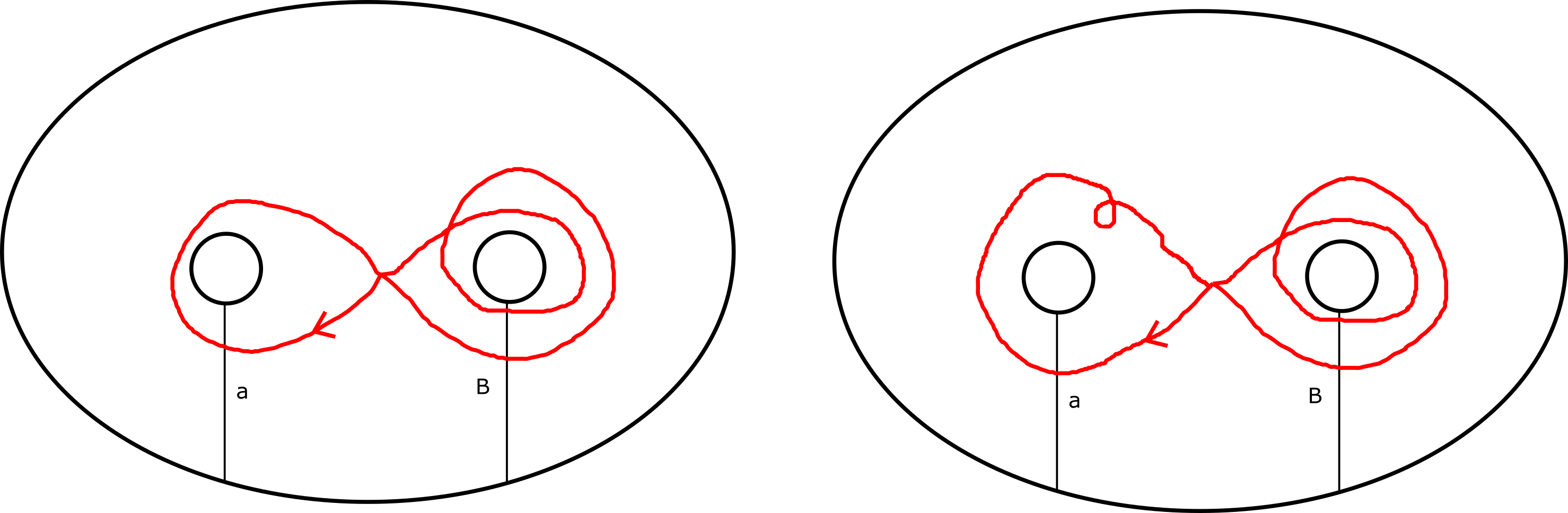}
    \caption{Both red closed curves are in the same deformation class.}
    \label{Basic}
\end{figure}

\begin{defn}\label{class} We say that two closed curves are in the same {\it deformation class} or {\it free homotopy class} if and only if the two closed curves can be continuously deformed into one another without leaving the surface. See Figure \ref{Basic} for an example.

\end{defn}

\begin{rem} In this paper, when we use the phrase {\it class of curves}, we always refer to the deformation class of curves or free homotopy class of curves defined in Definition \ref{class}, except when we clearly use the phrase {\it  equivalence class,} which is not defined until Section 3.3 with Definition \ref{triple}.

\end{rem}

\begin{rem} The word {\it curve} always refers to closed curves unless otherwise indicated.

\end{rem}

\begin{defn}\label{alphabet} Each class of curves can be described by a word. As in Figure \ref{Basic} above, we follow the curve based on its orientation and write down a, A, b or B depending on how the curve intersects the segments labeled a and B above. We always write the shortest form of the word, so aA and bB are never allowed. Further, we use C to denote ab and c to denote BA based on the relation abc = 1.

\end{defn}

\begin{defn} For two deformation classes of curves $\alpha, \beta$ we denote the minimum number of possible intersection points between a curve in class $\al$ and a curve in class $\beta$ by $\I(\al, \beta),$ counted with multiplicity. 

\end{defn}

\begin{defn} For any class of curves $c,$ we let $\si(c)$ be the minimal self intersection among all of the curves in $c.$ If $\si(c) = 0,$ we say that $c$ is a {\it simple closed curve}.

\end{defn}

\begin{defn}\label{equivalence} On a given surface $M$, two classes of curves $\al_1$ and $\al_2$ are defined to be {\it k-equivalent} if every class of curves $c$ on $M$ with self intersection number $k$ satisfies $i(\al_1, c) = i(\al_2, c).$

\end{defn}

\bexmp In Figure \ref{1eq} below, the three red curves each represent one of the classes of curves on the pair of pants with self intersection number $1,$ while we see that the intersection number between the blue and red curves is $2.$ 

\eexmp

\begin{figure}[H]
    \centering
    \includegraphics[scale = 0.2]{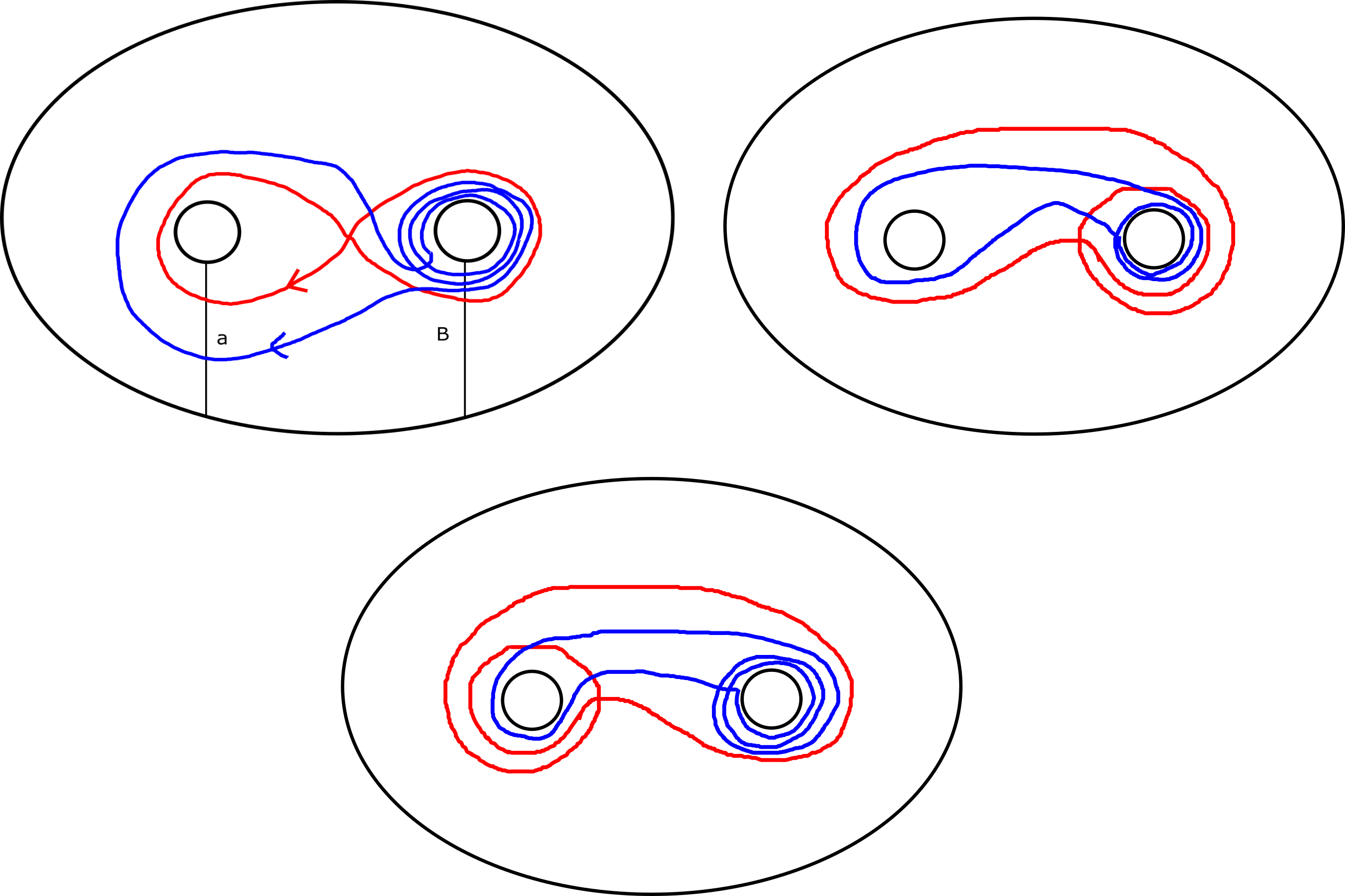}
    \caption{The blue curves represent $aB^n$ for integers $n \geq 1$ and are all 1-equivalent.}
    \label{1eq}
\end{figure}

Definition \ref{equivalence} is a generalization of the simple-intersection equivalence defined in \cite{leininger}.

\begin{defn}\label{power} A class $\gamma$ of curves is a {\it power} if there exists a class of curves $d$ and an integer $n \geq 2$ such that $d^n \in \gamma.$ Otherwise, $\gamma$ is a {\it non-power}.

\end{defn}

\begin{figure}[H]
    \centering
    \includegraphics[scale = 0.2]{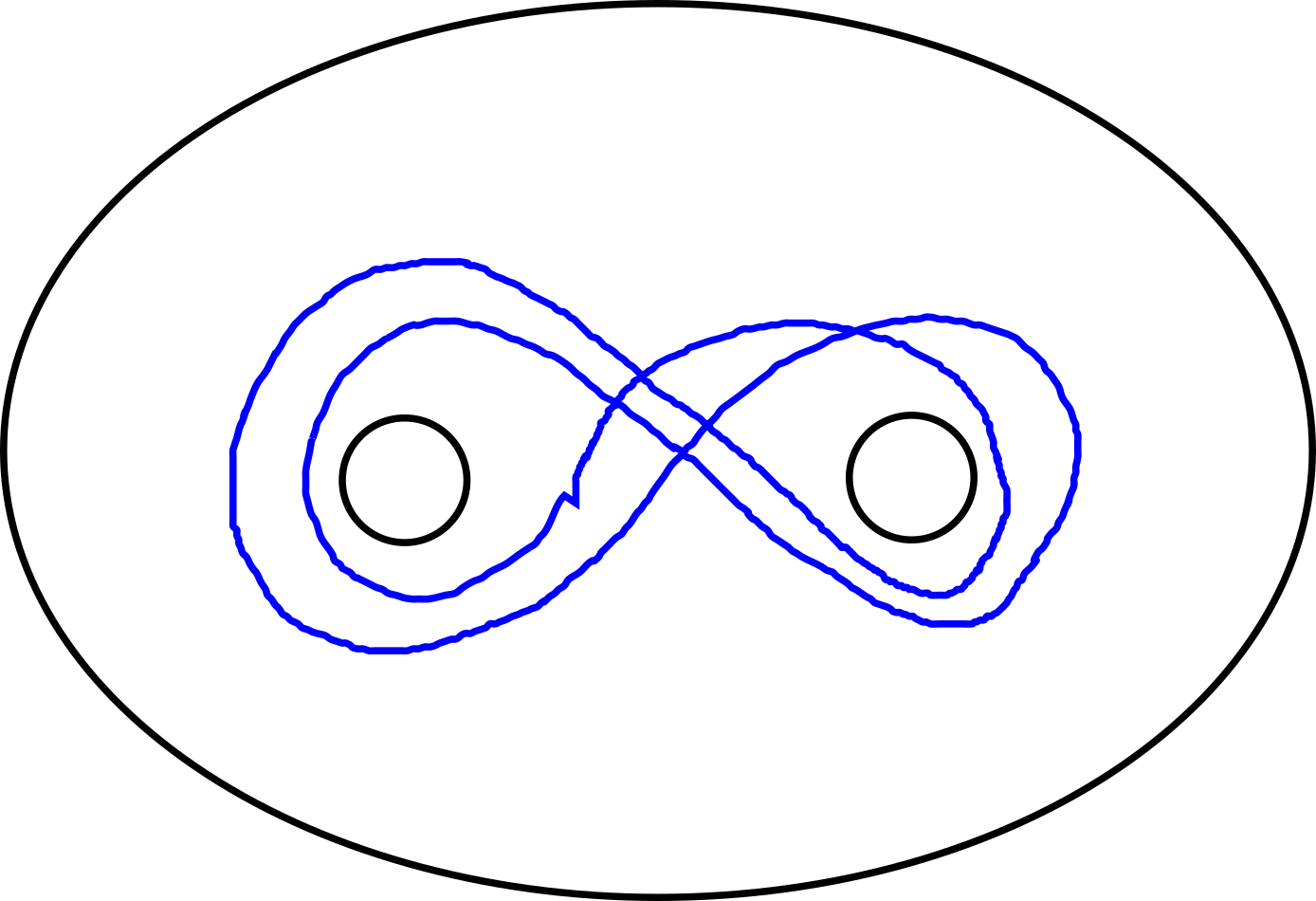}
    \caption{The curve aBaB, a power.}
    \label{powerfig}
\end{figure}

\begin{exmp} The blue curve in Figure \ref{powerfig} represents the second power of the figure $8$ curve aB on the pair of pants.

\end{exmp}

When we refer to classes of closed curves, we assume they are non-powers unless otherwise stated,

We will focus exclusively on powers of curves in the following section, but after that we will refer only to non-powers for the reasons stated in convention \ref{nopow}.




\ifx\allfiles\undefined

\newpage

\bibliography{Index}

\end{document}

\fi


\ifx\allfiles\undefined

\documentclass[12pt,a4paper]{article}


\usepackage{graphicx}

\usepackage{amsmath}

\usepackage{amssymb}

\usepackage{amsthm}

\usepackage{geometry}

\usepackage{fancyhdr}

\usepackage{color} 









\newcommand{\al}{\alpha}

\newcommand{\vphi}{\varphi}

\newcommand{\be}{\beta}

\newcommand{\ga}{\gamma}

\newcommand{\de}{\delta}

\newcommand{\om}{\omega}

\newcommand{\na}{\nabla}

\newcommand{\NA}{\nabla}

\newcommand{\bs}{\boldsymbol}

\newcommand{\ra}{\rightarrow}

\newcommand{\lra}{\longrightarrow}

\newcommand{\Ra}{\Rightarrow}

\newcommand{\xra}{\xrightarrow}

\newcommand{\xlra}{\xlongrightarrow}

\newcommand{\rgl}{\rangle}

\newcommand{\lgl}{\langle}

\newcommand{\dash}{\textrm{-}}

\newcommand{\ot}{\otimes}

\newcommand{\bpf}{\begin{proof}}

\newcommand{\epf}{\end{proof}}

\newcommand{\bthm}{\begin{thm}}

\newcommand{\ethm}{\end{thm}}

\newcommand{\bprop}{\begin{prop}}

\newcommand{\eprop}{\end{prop}}

\newcommand{\bcor}{\begin{cor}}

\newcommand{\ecor}{\end{cor}}

\newcommand{\blem}{\begin{lem}}

\newcommand{\elem}{\end{lem}}

\newcommand{\bdefn}{\begin{defn}}

\newcommand{\edefn}{\end{defn}}

\newcommand{\bexmp}{\begin{exmp}}

\newcommand{\eexmp}{\end{exmp}}

\newcommand{\brem}{\begin{rem}}

\newcommand{\erem}{\end{rem}}

\newcommand{\bdia}{\begin{displaymath}\xymatrix}

\newcommand{\edia}{\end{displaymath}}

\newcommand{\beq}{\begin{equation*}\begin{aligned}}

\newcommand{\eeq}{\end{aligned}\end{equation*}}

\newcommand{\bref}{\textbf{Ref}}

\newcommand{\intg}{\mathbb{Z}}

\newcommand{\real}{\mathbb{R}}

\newcommand{\comp}{\mathbb{C}}

\newcommand{\quot}{\mathbb{H}}

\DeclareMathOperator{\tr}{trunk}
\newcommand{\afv}{\mathbb{A}}

\newcommand{\prv}{\mathbb{P}}

\newcommand{\mco}{\mathcal{O}}

\newcommand{\mcc}{\mathcal{C}}

\newcommand{\mcf}{\mathcal{F}}

\newcommand{\mcg}{\mathcal{G}}

\newcommand{\mcs}{\mathcal{S}}

\newcommand{\cp}{\mathbb{CP}}

\newcommand{\mfo}{\mathfrak{O}}

\newcommand{\mfg}{\mathfrak{g}}

\newcommand{\msa}{\mathscr{A}}

\newcommand{\msr}{\mathscr{R}}

\newcommand{\msg}{\mathscr{G}}

\newcommand{\msd}{\mathscr{D}}

\newcommand{\itbf}{\item\textbf}

\newcommand{\seqa}{a_1,...,a_}

\newcommand{\seqx}{x_1,...,x_}

\newcommand{\seqy}{y_1,...,y_}

\newcommand{\seqf}{f_1,...,f_}

\newcommand{\cred}{\textcolor{red}}

\newcommand{\cblue}{\textcolor{blue}}

\newcommand{\mfa}{\mathfrak{a}}

\newcommand{\mfb}{\mathfrak{b}}

\newcommand{\mfm}{\mathfrak{m}}

\newcommand{\mfn}{\mathfrak{n}}

\newcommand{\mfp}{\mathfrak{p}}

\newcommand{\Af}{A_{(f)}}


\newtheorem{thm}{\textbf {Theorem}}[section]

\newtheorem{cor}[thm]{\textbf{Corollary}}

\newtheorem{prop}[thm]{\textbf{Proposition}}

\newtheorem{lem}[thm]{\textbf{Lemma}}

\newtheorem{conj}[thm]{Conjecture}

\newtheorem{prob}[thm]{Problem}

\newtheorem{exer}[thm]{Exercise}

\newtheorem{quest}[thm]{Question}

\theoremstyle{definition}

\newtheorem{defn}[thm]{\textbf{Definition}}

\newtheorem{defns}[thm]{Definitions}

\newtheorem{exmp}[thm]{Example}

\newtheorem{exmps}[thm]{Examples}

\newtheorem{var}[thm]{Variant}

\newtheorem{vars}[thm]{Variants}

\newtheorem{con}[thm]{Construction}

\newtheorem{notn}[thm]{Notation}

\newtheorem{notns}[thm]{Notations}

\newtheorem{conv}[thm]{Convention}

\theoremstyle{remark}

\newtheorem{rem}[thm]{Remark}

\newtheorem{rems}[thm]{Remarks}

\newtheorem{warn}[thm]{Warning}

\newtheorem{sch}[thm]{Scholium}

\newtheorem{expl}[thm]{Explanations}

\newtheorem*{theorem}{\textbf{Theorem}}

\newtheorem*{corollary}{\textbf{Corollary}}

\newtheorem*{proposition}{\textbf{Proposition}}

\newtheorem*{lemma}{\textbf{Lemma}}

\newtheorem*{example}{\textbf{Example}}

\def\cok{\operatorname{Coker}}

\newcommand{\txi}{\tilde{\xi}}

\newcommand{\bxi}{\bar{\xi}}

\newcommand{\bz}{\bar{z}}

\DeclareMathOperator{\tr}{tr}



\begin{document}

\bibliographystyle{plain}

\else

\fi

\section{A Study of k-Equivalence on the Pair of Pants}

\subsection{Powers of Curves and k-equivalence}

In this section, we examine how the properties of k-equivalence that were defined in definition \ref{equivalence} apply to the powers of curves defined in definition \ref{power}.

\begin{defn}\label{seq} Define the function $a(\ell, n) = ln^2 + n - 1.$

\end{defn}

\begin{lem}\label{seqlem} For a class of curves $\al$ such that $\si(\al) = \ell,$ it follows that $\si(\al^n) = a(\ell, n).$

\end{lem}

\begin{lem}\label{intpower} For classes of curves $\al, \beta$,  $\I(\al, \beta) = m$ if and only if $\I(\al^n, \beta) = mn.$

\end{lem}

\begin{prop}\label{powerprop} Suppose two classes of curves $\al_1, \al_2$ are $a(\ell, n)-$ equivalent for some positive integer $n$. Then it follows that $\al_1, \al_2$ are also $a(\ell, 1) = \ell$ equivalent.

\end{prop}

\bpf Consider a class of curves $\beta$ such that $\si(\beta) = \ell.$ Then, we know from lemma \ref{seqlem} that $\si(\beta^n) = a(\ell, n).$

Since $\al_1, \al_2$ are assumed to be $a(\ell, n)$ equivalent, we know that $\I(\al_1, \beta^n) = \I(\al_2, \beta^n).$ From lemma \ref{intpower}, it follows that $\I(\al_1, \beta) = \I(\al_2, \beta).$ Since $\beta$ was an arbitrary curve of self intersection number $\ell,$ it follows that $\al_1, \al_2$ must be $\ell$-equivalent as desired.

\epf

\subsection{Implications of Classes of Curves Being k-Equivalent}

In this section we use properties of classes of curves to prove the following equivalence relation:

\textbf{Theorem.} On the pair of pants, if two classes of curves $\al_1$ and $\al_2$ are k-equivalent for some positive integer $k \geq 2,$ then $\al_1$ and $\al_2$ are 2-equivalent and 1-equivalent. 

We begin by proving that if two classes of curves are k-equivalent then they must also be 1-equivalent. We then use k-equivalence and 1-equivalence together to prove 2-equivalence. We also provide evidence to suggest a similar argument can be used to prove that k-equivalence implies (k - 1)-equivalence on the pair of pants.

\begin{defn} On a fixed depiction of the pair of pants, we let $a^n$ denote a representative of the corresponding power formed by revolving $n$ times around component $a$ such that the two endpoints of the curve are collinear with the center of the boundary component, and then joining those endpoints to get a closed curve. We define $b^n$ and $c^n$ analogously.

\end{defn}



\begin{figure}[H]
    \centering
    \includegraphics[scale = .2]{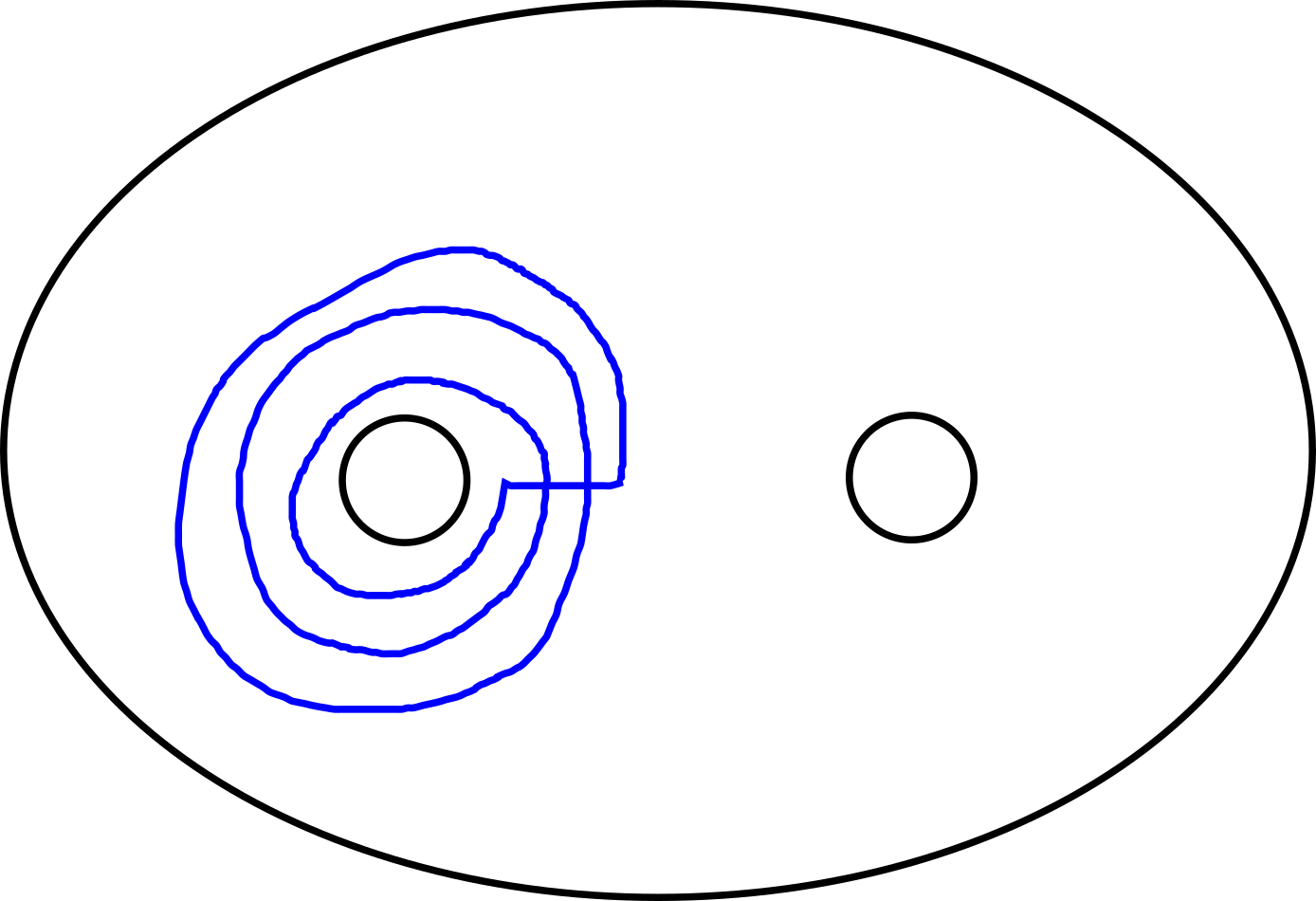}
    \caption{An example of $a^n$ on the pair of pants surface.}
    \label{xnpic}
\end{figure}

\begin{lem}[Chas, McMullen, Phillips \cite{moira}]\label{1fact} The only classes of closed curves on the pair of pants with one self intersection are aB, aC, Cb.

\end{lem}

\begin{lem}\label{self} For any positive integer $n,$ we have $\si(a^nB) = n, \si(a^nCb) = n + 1$ and $\si(a^nCC) = n + 1.$

\end{lem}


\begin{figure}[H]
    \centering
    \includegraphics[scale = 0.125]{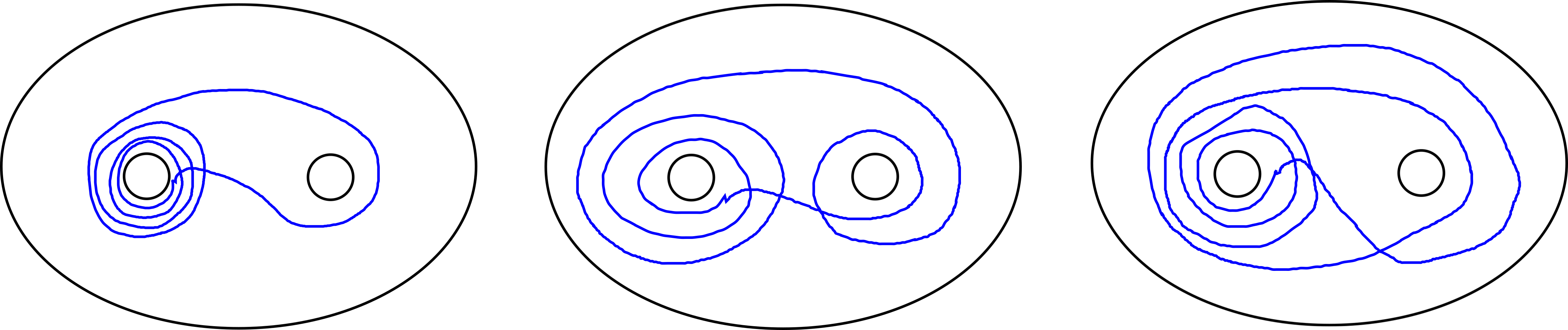}
    \caption{$a^nB, a^nCb, a^nCC$}
    \label{fig:my_label}
\end{figure}



\begin{notn} We use the notation $\Int(m_1, m_2)$ to denote the intersection number of two representatives $m_1, m_2$ of classes of curves.  

\end{notn}

\begin{lem}\label{kimp11} Suppose on a pair of pants, we have two classes of curves $\al_1$ and $\al_2$ such that the following holds:

$\I(\al_1, a^nB) = \I(\al_2, a^nB)$

$\I(\al_1, a^nC) = \I(\al_2, a^nC)$

$\I(\al_1, a^nbC) = \I(\al_2, a^nbC)$

Then it follows that $\I(\al_1, Cb) = \I(\al_2, Cb).$
\end{lem}

\bpf 

Let $m_1, m_2$ be respective representatives of classes $\al_1, \al_2$ such that $\Int(m_1, a^nB)$, $\Int(m_1, a^nC)$, $\Int(m_1, a^nbC)$, $\Int(m_2, a^nB)$, $\Int(m_2, a^nC)$, $\Int(m_2, a^nbC)$ are minimized.

We observe that there exists a set of representatives of $a^nB, a^nC, a^nbC$ with minimal self intersection such that the $a^n$ part of the curves coincide, and the loops around the $b$ and $c$ boundary components coincide as well. 

Since we are dealing with representatives of classes of curves which coincide as described above, it follows that $m_1, m_2$ will have a fixed number of intersections with the representatives of $a^n$ and the loops around the $b$ and $c$ components.

Let $a_1, b_1, c_1$ be the intersection numbers between the curves $a^n, b, C$ and $m_1.$ Define $a_2, b_2, c_2$ analogously. By translating the assumed intersection relations between classes of curves above to representatives of those classes of curves, we get the following system of equations:

$$a_1 + b_1 = a_2 + b_2$$ 
$$a_1 + c_1 = a_2 + c_2$$
$$a_1 + b_1 + c_1 = a_2 + b_2 + c_2.$$

From this, it clearly follows that $a_1 = a_2, b_1 = b_2, c_1 = c_2.$ Therefore, we get that $\I(\al_1, Cb) = \I(\al_2, Cb).$ 

\epf

\begin{rem}\label{eqint} Note that the solution $a_1 = a_2, b_1 = b_2,c_1=c_2$ means that the corresponding intersection numbers between each of the loops and the deformation classes $\al_1, \al_2$ are the same. We will use this fact to complete the proof of the theorem at the beginning of this subsection.

\end{rem}

\begin{rem} The exact same method used in the proof of lemma \ref{kimp11} can be also used to prove the following two lemmas:

\end{rem}

\begin{lem}\label{kimp12} Suppose on a pair of pants, we have two classes of curves $\al_1$ and $\al_2$ such that the following holds:

$\I(\al_1, b^nC) = \I(\al_2, b^nC)$

$\I(\al_1, b^nA) = \I(\al_2, b^nA)$

$\I(\al_1, b^naC) = \I(\al_2, b^naC)$

Then it follows that $\I(\al_1, aC) = \I(\al_2, aC).$
\end{lem}

\begin{lem}\label{kimp13} Suppose on a pair of pants, we have two classes of curves $\al_1$ and $\al_2$ such that the following holds:

$\I(\al_1, c^nB) = \I(\al_2, c^nB)$

$\I(\al_1, c^nA) = \I(\al_2, c^nA)$

$\I(\al_1, c^naB) = \I(\al_2, c^naB)$

Then it follows that $\I(\al_1, aB) = \I(\al_2, aB).$
\end{lem}

\begin{thm}\label{kimp1thm} On a pair of pants, if two classes of curves $\al_1, \al_2$ are $k$-equivalent for some integer $k \geq 2,$ it follows that $\al_1$ and $\al_2$ are also $1$-equivalent.

\end{thm}

\bpf This follows from lemmas \ref{kimp11}, \ref{kimp12} and \ref{kimp13}.

\epf 

\begin{lem} The only deformation classes of curves on the pair of pants with self intersection number $2$ are aaB, aBB, aaC, aCC, Cbb, CCb, aCb, CaB, CAb.

\end{lem}

\begin{prop}\label{kimp21} If two classes of curves $\al_1, \al_2$ on a pair of pants are k-equivalent for some positive integer $k \geq 2,$ we get the following:

$\I(\al_1, aBB) = \I(\al_2, aBB).$

\begin{figure}[H]
    \centering
    \includegraphics[scale = 0.2]{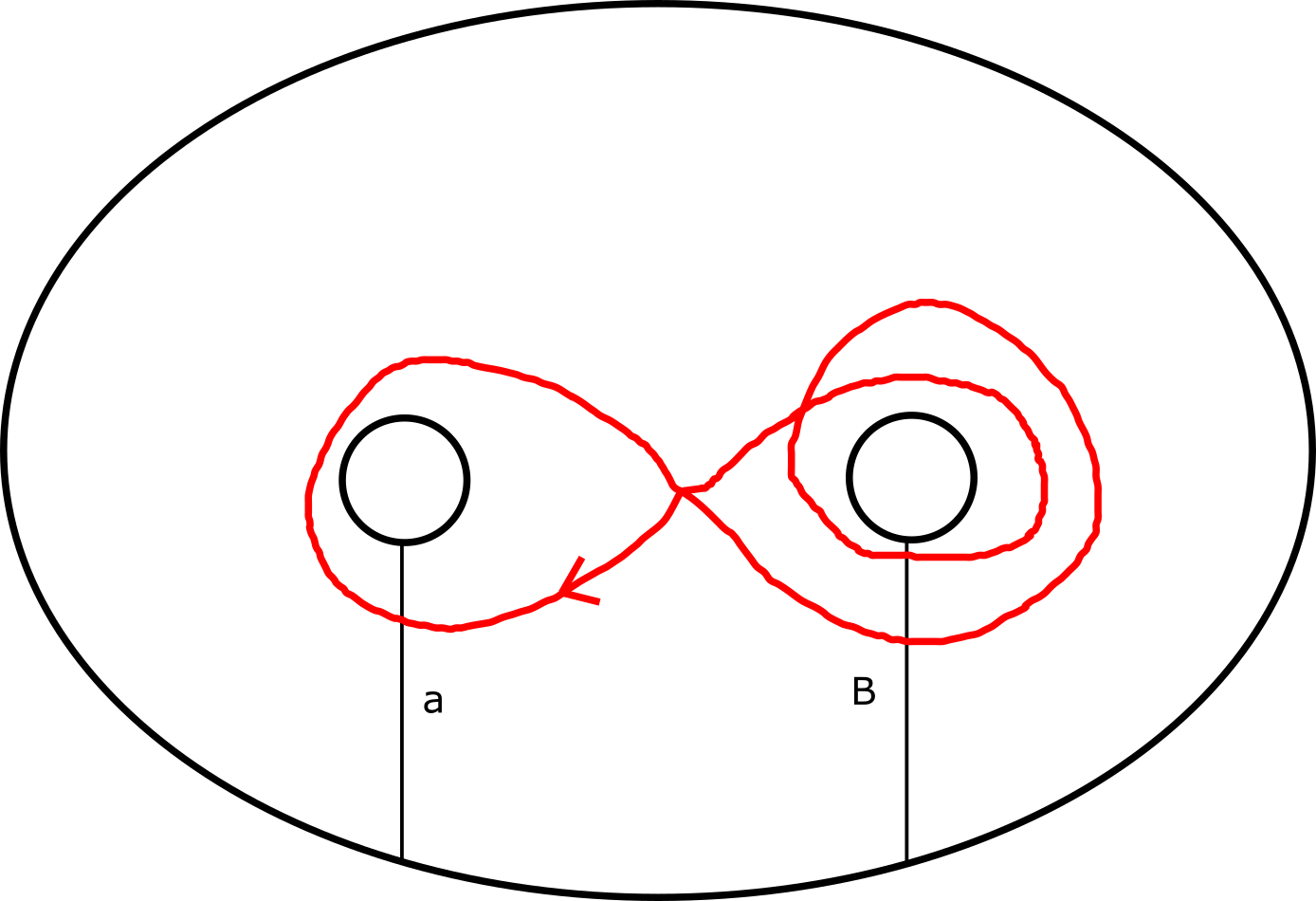}
    \caption{A representative of the class of curves aBB on the pair of pants.}
    \label{aBB}
\end{figure}

\bpf


From remark \ref{eqint}, we know that both $\al_1$ and $\al_2$ will intersect the loop around the ``a" boundary component the same number of times, and they will also each intersect $B^2$ the same number of times, so we are done.



\epf

\end{prop}

\begin{rem}\label{2eqint} The proof of proposition \ref{kimp21} can be also used to show the same result for aaC, aCC, aaB, CCb and Cbb, all of which have self intersection number $2.$ Each of these classes of curves has the same topological nature as aBB on the pair of pants. 

\end{rem}

Using theorem \ref{kimp1thm}, we can now prove the following:

\begin{thm}\label{kimp2thm} On a pair of pants, if two classes of curves $\al_1, \al_2$ are k-equivalent, it follows that they are also 2-equivalent.

\end{thm}

\bpf

It remains to only prove the equivalence relation for aCb, CaB and CAb, as all other classes of curves with self intersection number $2$ have been addressed by remark \ref{2eqint}. We begin with aCb.

\begin{figure}[H]
    \centering
    \includegraphics[scale = 0.2]{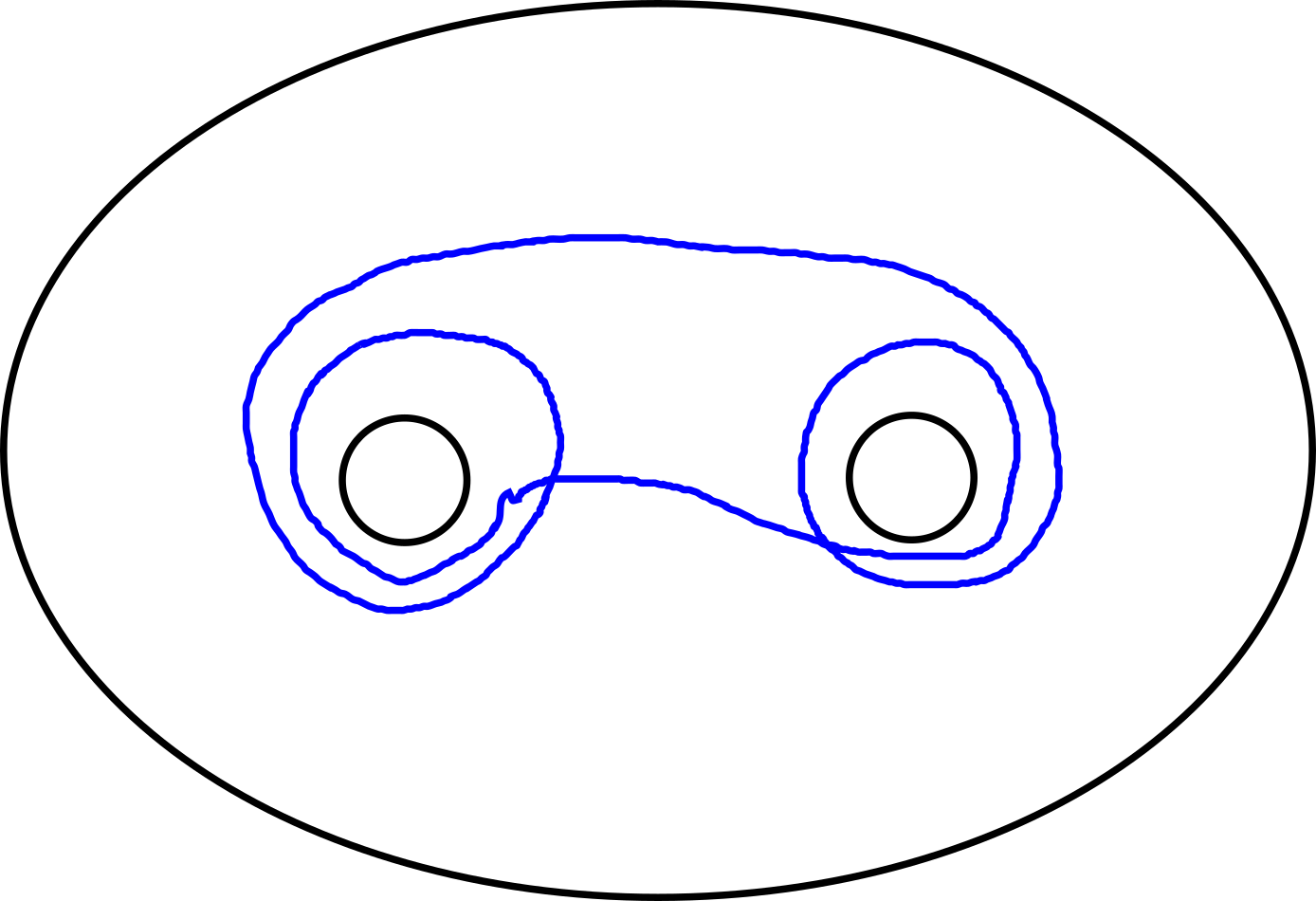}
    \caption{A representative of the class of curves aCb on the pair of pants.}
    \label{aCb}
\end{figure}

We know that the $\al_1, \al_2$ are 1-equivalent from theorem \ref{kimp1thm}. We fix representatives of $\al_1, \al_2, aCb$ so that the intersection numbers between the representatives of $\al_1, \al_2$ with $aCb$ are minimized. Then, from 1-equivalence, we know that the representatives will each intersect the figure 8 curve Cb the same number of times. Further, we know that the two classes of curves intersect the a-loop the same number of times. It follows that $\I(\al_1, aCb) = \I(\al_2, aCb).$

Representatives of the classes of curves CaB and CAb are below:

\begin{figure}[H]
    \centering
    \includegraphics[scale = 0.2]{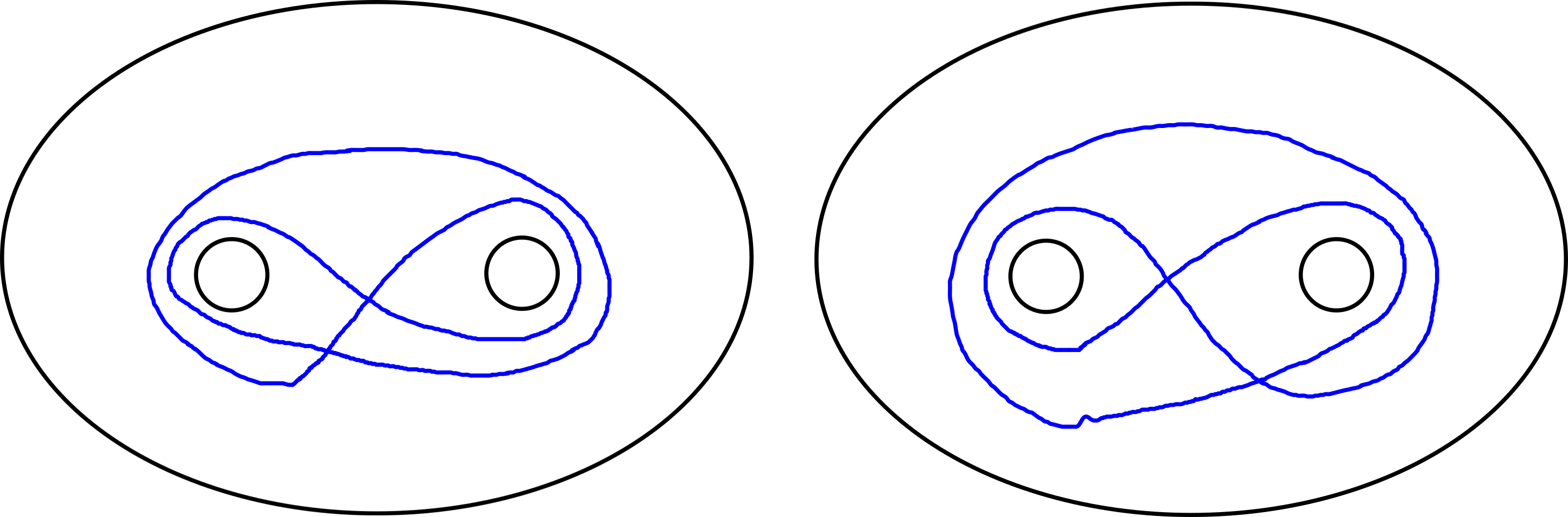}
    \caption{The left curve is a representative of CaB, while the right curve is a representative of CAb.}
    \label{2si}
\end{figure}

Both of these representatives of classes of curves can be viewed as a combination of the aB figure $8$ and the C loop. As $\al_1, \al_2$ are 1-equivalent, we get that $\I(\al_1, CAb) = \I(\al_2, CAb)$ and $\I(\al_1, CaB) = \I(\al_2, CaB)$ by the same argument as above. 

\epf

\subsection{Patterns in 1-Equivalence}

In this section we examine a phenomenon of 1-equivalence. Recall that there are exactly three curves on the pair of pants with self intersection number $1:$ aB, Cb, aC.

\begin{defn}\label{triple} For any curve $\al$ on the pair of pants, we have an associated ordered triple $t(\al)= \left(\I(\al, aB), \I(\al, Cb), \I(\al, aC)\right).$ We call $t(\al)$ an {\it equivalence class}, as it denotes the set of classes of curves that are 1-equivalent to $\al$ on the pair of pants.

\end{defn}

\begin{lem}
Suppose a class of curves $l$ has two intersections with $aB$ on a pair of pants. Then $l$ must be one of the following forms:

$$C^ma^nb, C^ma^nB, C^mab^n, C^maB^n, C^mA^nb, C^mA^nB, C^mAb^n, C^mAB^n$$ for some positive integer $n.$

\end{lem}

\bpf

A representative $d$ of the class $aB$ has a loop around each of the ``small" boundary components as shown in Figure \ref{aB} below:

\begin{figure}[H]
    \centering
    \includegraphics[scale = 0.2]{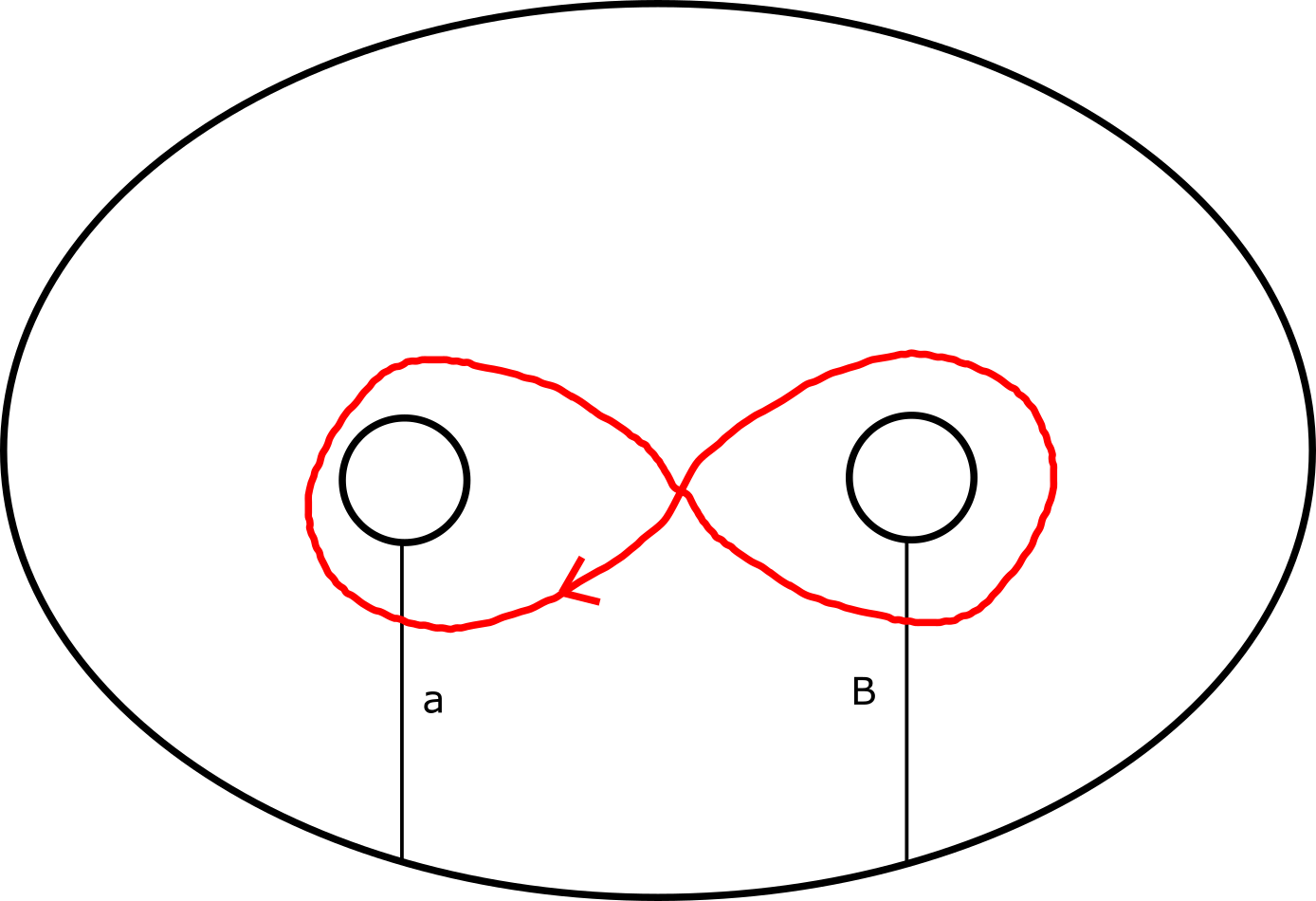}
    \caption{A representative of one of the classes of figure 8 curves on the pair of pants.}
    \label{aB}
\end{figure}

Suppose arbitrarily that $l$ intersects the loop around the leftmost boundary component (the proof is analogous if $l$ intersects the other loop). Then, it can either go around the $a$ boundary component clockwise or counterclockwise any number of times before exiting. Once it exits, the curve can loop around $ab$ any number of times and then eventually come back to $a.$

\epf 

\begin{figure}[H]
    \centering
    \includegraphics[scale = 0.2]{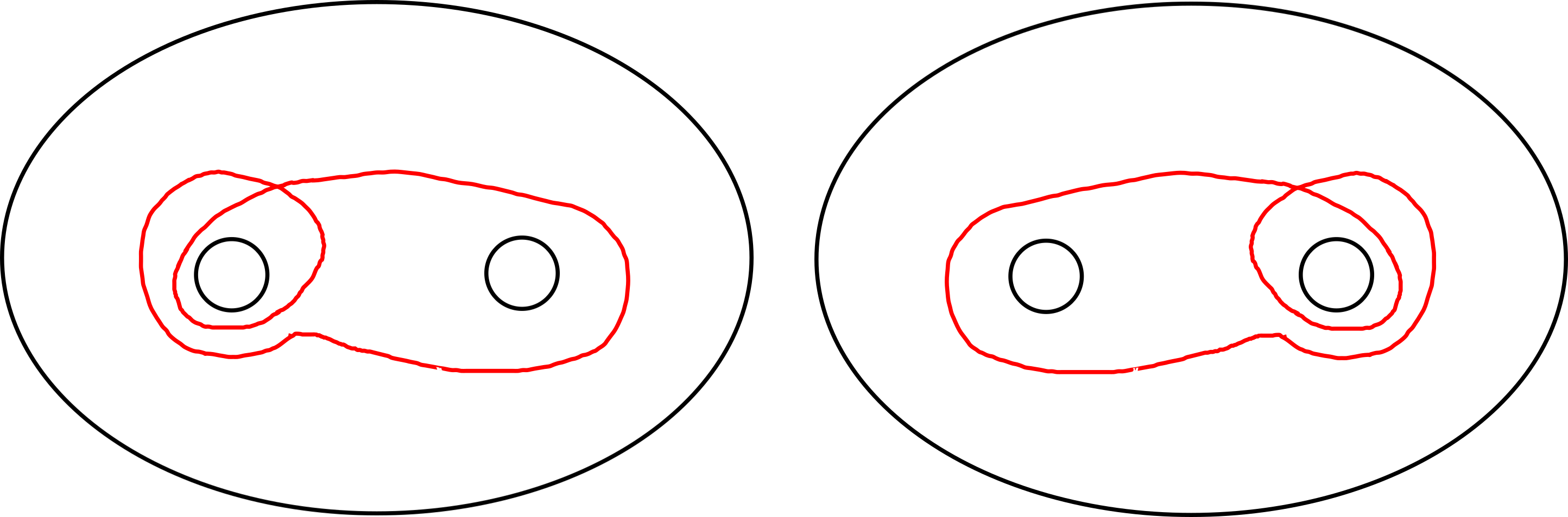}
    \caption{Representatives of the other two classes of curves with self intersection number $1$: Cb and aC}
    \label{CbaC}
\end{figure}

\begin{prop}\label{classification} If a class of curves $\al$ is in the equivalence class $t(\al) = (2, 2, 2),$ then $\al$ must be one of the following forms: 

\end{prop}

\begin{conj}\label{1equivdouble} For any curve $\al$ on the pair of pants, $\max(t(\al)) \leq 2 \cdot \min(t(\al)).$ Further, $\max(t(\al)) = 2 \cdot \min(t(\al))$ if and only if $t(\al)$ is some permutation of $(q, q, 2q)$ where $q$ is an even positive integer. We require $q$ to be even because of the fact that any the intersection number of any two deformation classes of curves on the pair of pants must be even.

\end{conj}

Computer experimentation reveals that $t(\al)$ can equal (2, 2, 2,), (2, 2, 4), (4, 4, 4), (4, 4, 6), (4, 4, 8) and so on, including all permutations of those listed.




\ifx\allfiles\undefined

\newpage

\bibliography{Index}

\end{document}

\fi


\ifx\allfiles\undefined

\documentclass[12pt,a4paper]{article}


\usepackage{graphicx}

\usepackage{amsmath}

\usepackage{amssymb}

\usepackage{amsthm}

\usepackage{geometry}

\usepackage{fancyhdr}

\usepackage{color} 









\newcommand{\al}{\alpha}

\newcommand{\vphi}{\varphi}

\newcommand{\be}{\beta}

\newcommand{\ga}{\gamma}

\newcommand{\de}{\delta}

\newcommand{\om}{\omega}

\newcommand{\na}{\nabla}

\newcommand{\NA}{\nabla}

\newcommand{\bs}{\boldsymbol}

\newcommand{\ra}{\rightarrow}

\newcommand{\lra}{\longrightarrow}

\newcommand{\Ra}{\Rightarrow}

\newcommand{\xra}{\xrightarrow}

\newcommand{\xlra}{\xlongrightarrow}

\newcommand{\rgl}{\rangle}

\newcommand{\lgl}{\langle}

\newcommand{\dash}{\textrm{-}}

\newcommand{\ot}{\otimes}

\newcommand{\bpf}{\begin{proof}}

\newcommand{\epf}{\end{proof}}

\newcommand{\bthm}{\begin{thm}}

\newcommand{\ethm}{\end{thm}}

\newcommand{\bprop}{\begin{prop}}

\newcommand{\eprop}{\end{prop}}

\newcommand{\bcor}{\begin{cor}}

\newcommand{\ecor}{\end{cor}}

\newcommand{\blem}{\begin{lem}}

\newcommand{\elem}{\end{lem}}

\newcommand{\bdefn}{\begin{defn}}

\newcommand{\edefn}{\end{defn}}

\newcommand{\bexmp}{\begin{exmp}}

\newcommand{\eexmp}{\end{exmp}}

\newcommand{\brem}{\begin{rem}}

\newcommand{\erem}{\end{rem}}

\newcommand{\bdia}{\begin{displaymath}\xymatrix}

\newcommand{\edia}{\end{displaymath}}

\newcommand{\beq}{\begin{equation*}\begin{aligned}}

\newcommand{\eeq}{\end{aligned}\end{equation*}}

\newcommand{\bref}{\textbf{Ref}}

\newcommand{\intg}{\mathbb{Z}}

\newcommand{\real}{\mathbb{R}}

\newcommand{\comp}{\mathbb{C}}

\newcommand{\quot}{\mathbb{H}}

\DeclareMathOperator{\tr}{trunk}
\newcommand{\afv}{\mathbb{A}}

\newcommand{\prv}{\mathbb{P}}

\newcommand{\mco}{\mathcal{O}}

\newcommand{\mcc}{\mathcal{C}}

\newcommand{\mcf}{\mathcal{F}}

\newcommand{\mcg}{\mathcal{G}}

\newcommand{\mcs}{\mathcal{S}}

\newcommand{\cp}{\mathbb{CP}}

\newcommand{\mfo}{\mathfrak{O}}

\newcommand{\mfg}{\mathfrak{g}}

\newcommand{\msa}{\mathscr{A}}

\newcommand{\msr}{\mathscr{R}}

\newcommand{\msg}{\mathscr{G}}

\newcommand{\msd}{\mathscr{D}}

\newcommand{\itbf}{\item\textbf}

\newcommand{\seqa}{a_1,...,a_}

\newcommand{\seqx}{x_1,...,x_}

\newcommand{\seqy}{y_1,...,y_}

\newcommand{\seqf}{f_1,...,f_}

\newcommand{\cred}{\textcolor{red}}

\newcommand{\cblue}{\textcolor{blue}}

\newcommand{\mfa}{\mathfrak{a}}

\newcommand{\mfb}{\mathfrak{b}}

\newcommand{\mfm}{\mathfrak{m}}

\newcommand{\mfn}{\mathfrak{n}}

\newcommand{\mfp}{\mathfrak{p}}

\newcommand{\Af}{A_{(f)}}


\newtheorem{thm}{\textbf {Theorem}}[section]

\newtheorem{cor}[thm]{\textbf{Corollary}}

\newtheorem{prop}[thm]{\textbf{Proposition}}

\newtheorem{lem}[thm]{\textbf{Lemma}}

\newtheorem{conj}[thm]{Conjecture}

\newtheorem{prob}[thm]{Problem}

\newtheorem{exer}[thm]{Exercise}

\newtheorem{quest}[thm]{Question}

\theoremstyle{definition}

\newtheorem{defn}[thm]{\textbf{Definition}}

\newtheorem{defns}[thm]{Definitions}

\newtheorem{exmp}[thm]{Example}

\newtheorem{exmps}[thm]{Examples}

\newtheorem{var}[thm]{Variant}

\newtheorem{vars}[thm]{Variants}

\newtheorem{con}[thm]{Construction}

\newtheorem{notn}[thm]{Notation}

\newtheorem{notns}[thm]{Notations}

\newtheorem{conv}[thm]{Convention}

\theoremstyle{remark}

\newtheorem{rem}[thm]{Remark}

\newtheorem{rems}[thm]{Remarks}

\newtheorem{warn}[thm]{Warning}

\newtheorem{sch}[thm]{Scholium}

\newtheorem{expl}[thm]{Explanations}

\newtheorem*{theorem}{\textbf{Theorem}}

\newtheorem*{corollary}{\textbf{Corollary}}

\newtheorem*{proposition}{\textbf{Proposition}}

\newtheorem*{lemma}{\textbf{Lemma}}

\newtheorem*{example}{\textbf{Example}}

\def\cok{\operatorname{Coker}}

\newcommand{\txi}{\tilde{\xi}}

\newcommand{\bxi}{\bar{\xi}}

\newcommand{\bz}{\bar{z}}

\DeclareMathOperator{\tr}{tr}



\begin{document}

\bibliographystyle{plain}

\else

\fi

\section{Conclusion}

\begin{subsection}{Summary}

With proposition \ref{powerprop}, we see how k-equivalence works when powers of classes of curves are considered, as the self intersection number of a power and the intersection number of a power with another class of curves are both well understood. After that, we only look at non-powers of curves. With theorem \ref{kimp1thm}, we first prove that k-equivalence implies 1-equivalence on the pair of pants. Then, with theorem \ref{kimp2thm}, we extend this to show that 2-equivalence is implied as well. In Section 3.3, we fully describe the equivalence class $(2, 2, 2)$ with proposition \ref{classification}. We also see a pattern in 1-equivalence and conjecture that it holds in general for all classes of curves on the pair of pants. 

\end{subsection}

\begin{subsection}{Future Directions of Study}

The strongest conjecture regarding k-equivalence on the pair of pants would be the following:

\begin{conj} On the pair of pants, k-equivalence implies (k - 1)-equivalence for all positive integers $k.$

\end{conj}


Theorems \ref{kimp1thm} and \ref{kimp2thm} offer support towards this conjecture, along with the result that all classes of curves on the pair of pants are simple intersection equivalent. It appears possible to continue building up this logic to prove 3-equivalence, 4-equivalence all the way up to (k - 1)-equivalence. However, the number of classes of curves on a pair of pants with a fixed number of self intersections increases rapidly, making the methods used here difficult to generalize for large values of $k.$





We can also continue focusing on conjecture \ref{1equivdouble}, which gives us information about how different deformation classes of figure $8$ curves intersect other classes of curves. It is also possible that conjecture \ref{1equivdouble} has generalizations that extend beyond 1-equivalence.

It is also currently unknown if the equivalence relations proved here for the pair of pants surface are also true on other surfaces with boundary such as the punctured torus. However, other surfaces with boundary will be harder to study because only the pair of pants has a finite number of classes of curves for any fixed self intersection number.

\end{subsection}


\ifx\allfiles\undefined

\newpage

\bibliography{Index}

\end{document}

\fi



\end{document}